\documentclass[reqno,10pt]{amsart}
\usepackage{amsthm,amsmath,amsfonts,amssymb}
\usepackage{graphics,graphicx,epsf}
\usepackage{hyperref,color}         
\usepackage{dsfont}
\usepackage{epsfig}

\newcommand{\bproof}{\noindent\textbf{Proof: }}

\newcommand{\R}{{\mathbb R}}

\newcommand{\N}{{\mathbb N}}
\newcommand{\cqd}{\hspace*{\fill}\rule{0.2cm}{0.2cm}}  
\def\eps{\varepsilon}

\def\dist{{\rm dist}}
\def\eto{\buildrel \eps\to 0\over\longrightarrow }
\def\ds

\def \ts {\textstyle}
\def \ds {\displaystyle}

\def \demb {{{ \bigcap\kern -12.2pt{_\big\downarrow}}}}
\def \uemb {{{ \bigcup\kern -5.8pt{^\big\uparrow}}}}

\def \Re{{I \kern -.3em R}}
\def \RN {{{\bf I\kern -1.6pt{\bf R}}}^\N}
\def \Rn {{{\bf I\kern -1.6pt{\bf R}}}^{\rm n}}

\def \N {N}
\def \RN {\Re^N}

\def \ts {\textstyle}
\def \eps {\epsilon}
\def \la {\lambda}
\def \ga {\gamma}

\newtheorem{lemma}{Lemma}[section]

\newtheorem{proposition}[lemma]{Proposition}
\newtheorem{corollary}[lemma]{Corollary}

\newtheorem{definition}[lemma]{Definition}
\newtheorem{remark}[lemma]{Remark}

\numberwithin{equation}{section}

\begin{document}

\title[Geometric versus spectral convergence]%
{Geometric versus spectral convergence for the Neumann Laplacian
under exterior perturbations of the domain}

\author[J. M. Arrieta \& D. Krej\v{c}i\v{r}{\'\i}k]%
{Jos\'{e} M. Arrieta$^*$ and David Krej\v{c}i\v{r}{\'\i}k$^\dag$}
\thanks{$^*$Partially supported by Grants PHB2006-003-PC and  MTM2006-08262 from MEC and  by
``Programa de
Financiaci{\'o}n de Grupos de Investigaci{\'o}n UCM-Comunidad de Madrid
CCG07-UCM/ESP-2393. Grupo 920894'' and SIMUMAT-Comunidad de Madrid,  Spain
}

\thanks{$^\dag$Partially supported by the Czech Ministry of Education,
Youth and Sports within the project LC06002.
}

\address[J. M. Arrieta]
{Departamento de Matem{\'a}tica Aplicada, Universidad Complutense de
Madrid, 28040 Madrid  Spain.}
\email{arrieta@mat.ucm.es}

\address[D. Krej\v{c}i\v{r}{\'\i}k]
{Department of Theoretical Physics,
Nuclear Physics Institute,
Academy of Sciences,
25068 \v{R}e\v{z},
Czech Republic.}
\email{krejcirik@ujf.cas.cz}

\date{29 January 2009}

\begin{abstract}
We analyze the behavior of the eigenvalues and eigenfunctions of the
Laplace operator with homogeneous Neumann boundary conditions when
the domain is perturbed. We show that if $\Omega_0\subset\Omega_\eps$
are bounded domains (although not necessarily uniformly bounded) and
we know that the eigenvalues and eigenfunctions with Neumann boundary condition in
$\Omega_\eps$ converge to the ones in $\Omega_0$, then necessarily we have
that $|\Omega_\eps\setminus\Omega_0|\to 0$ while it is not necessarily true
that dist$(\Omega_\eps,\Omega_0)\eto 0$. As a matter of fact we will construct
an example of a perturbation where the spectra behave continuously but
 dist$(\Omega_\eps,\Omega_0)\eto +\infty$.
\end{abstract}

\maketitle

\section{Introduction}\label{Sintro}
This paper is concerned with the behavior of the eigenvalues
and eigenfunctions of the Laplace operator
in bounded domains when the domain undergoes a perturbation.
It is well known that if the boundary condition that
we are imposing is of Dirichlet type,
the kind of perturbations that we may allow in order to obtain
the continuity of the spectra is much broader than
in the case of Neumann boundary condition.
This is explicitly stated in the pioneer work of Courant and Hilbert \cite{CouHil37}
and it has been subsequently
clarified in many works, see \cite{BabVib65, Arr97, Dan03}
and reference therein among others.
See also \cite{Henrot2006} for a general text on
different properties of eigenvalues and \cite{Hen96}
for a study on the behavior of eigenvalues and in general partial
differential equations when the domain is perturbed.

In particular, with Dirichlet boundary condition
we may consider the case where the fixed domain is
a bounded ``smooth'' domain $\Omega_0\subset \R^N$, $N\geq 2$,
and the perturbed domain is $\Omega_\eps$ in such a
way that $\Omega_0\subset\Omega_\eps$,
that is we consider exterior perturbation of the domain. We
may have perturbations of this type where
$|\Omega_\eps\setminus\Omega_0|\geq \eta$ for some fixed
$\eta>0$ and still we have the convergence of the eigenvalues and eigenfunctions.
Moreover, we may even have the case
$|\Omega_\eps\setminus\Omega_0|\to +\infty$ and still we  have the convergence
of the eigenvalues and eigenfunctions.

To obtain and example of this situation is not too difficult.
If we consider for instance $\Omega\subset\R^2$,
given by  $\Omega_0=(0,1)\times (-1,0)$ and
$$\Omega_\eps(a)=\{(x,y):0<x<1, -1<y<a(1+\sin(x/\eps))\}\supset \Omega_0$$
where $a>0$ is fixed, we can easily see that the eigenvalues and eigenfunctions
of the Laplace operator
with Dirichlet boundary condition in $\Omega_\eps$ converge to the ones in $\Omega_0$.
Moreover
$
  |\Omega_\eps|=|\Omega_0|+ \int_0^1a(1+\sin(x/\eps))dx\sim |\Omega_0|+a
$
for $\eps$ small enough.  Moreover, it is not difficult to modify the example above
 choosing the constant $a$ dependent with respect to $\eps$
 in such a way that $a(\eps)\to +\infty$ and such that we still get
  that the eigenvalues and eigenfunctions in $\Omega_\eps(a(\eps))$ converge to the ones in $\Omega_0$ and
$|\Omega_\eps(a(\eps))\setminus\Omega_0|\to +\infty$.  This example shows that the class of perturbations that we
may allow to get the ``spectral convergence'' of the Dirichlet Laplacian is very broad and that knowing that the
eigenvalues and eigenfunctions of the Dirichlet Laplacian converge does not have many ``geometrical'' restrictions
for the domains.

The case of Neumann boundary condition is much more subtle. As a matter of fact, for the situation depicted above it is not
true that the spectra converge. So we ask ourselves the following questions: if we have a domain $\Omega_0$ and consider a perturbation of
it given by
$\Omega_0\subset\Omega_\eps$, where we assume that all the domains are smooth and bounded although not necessarily
uniformly bounded on the parameter $\eps$, then if we have the convergence of the eigenvalues and eigenfunctions,
\begin{itemize}

\item[{\bf (Q1)}] should it be true that $|\Omega_\eps\setminus\Omega_0|\eto 0$?

\item[]

\item[{\bf (Q2)}] should it be true that dist$(\Omega_\eps,\Omega_0)=\sup_{x\in\Omega_\eps}$dist$(x,\Omega_0)\eto 0$?

\end{itemize}

\par\bigskip

We will see that the answer to the first question is Yes  and, surprisingly,  the answer to the second one is No.

Observe that, as the example above shows, the answer to both questions for the case of Dirichlet boundary condition is
No.

In Section \ref{convergence} we recall a result from \cite{ar2,AC} which provides a necessary and sufficient condition for the convergence of eigenvalues
and eigenfunctions when the domain is perturbed. In Section \ref{Measure-convergence} we provide an answer to question {\bf(Q1)} and in Section \ref{distance-convergence} we provide an answer to question {\bf (Q2)}.

\section{Characterization of spectral convergence of Neumann Laplacian}
\label{convergence}

In this section we give a necessary and sufficient condition for the convergence of the eigenvalues and
eigenfunctions of the Laplace operator with Neumann boundary conditions. We refer to \cite{ar2} and \cite{AC} for
a general result in this direction, even in a more general context than the one in this note.
In our particular case, we will consider the following situation: let $\Omega_0$ be a fixed bounded
smooth (Lipschitz is enough) open set in $\R^N$ with $N\geq 2$ and let $\Omega_\eps$ be a family of domains such that
for each fixed $0<\eps\leq\eps_0$, $\Omega_\eps$ is bounded and smooth with
$ \Omega_0\subset\Omega_\eps$.

Let us define now what we mean by the spectral convergence.
For $0\leq\eps\leq\eps_0$, we
denote by $\{\lambda_n^\eps\}_{n=1}^\infty$
the sequence of eigenvalues of the Neumann Laplacian in $\Omega_\eps$,
always ordered and counting its multiplicity,
and we denote by $\{\phi_n^\eps\}_{n=1}^\infty$ a corresponding
set of orthonormal eigenfunctions in $\Omega_\eps$.
Also, since we are considering domains which vary with the parameter $\eps$
and we will need to compare functions defined in $\Omega_0$ and in $\Omega_\eps$,
we introduce the following space
$H^1_\eps=H^1(\Omega_0)\oplus H^1(\Omega_\eps\setminus\bar\Omega_0)$,
that is $\chi\in H^1_\eps$ if $\chi_{|\Omega_0}\in H^1(\Omega_0)$ and
$\chi_{|(\Omega_\eps\setminus\bar\Omega_0)}\in H^1(\Omega_\eps\setminus\bar\Omega_0)$,
with the norm
$$
  \|\chi\|_{H^1_\eps}^2
  = \|\chi\|_{H^1(\Omega_0)}^2+ \|\chi\|_{H^1(\Omega_\eps\setminus \bar\Omega_0)}^2
  \,.
$$

We have that $H^1(\Omega_\eps)\hookrightarrow H^1_\eps$
and in a natural way we have that if $\chi\in H^1(\Omega_0)$
via the extension by zero outside $\Omega_0$ we have $\chi\in H^1_\eps$.
Hence, with certain abuse of notation
we may say that if $\chi_\eps\in H^1_\eps$, $0\leq\eps\leq\eps_0$, then
$\chi_\eps\eto \chi_0$  in $H^1_\eps$ if
$\|\chi_\eps-\chi_0\|_{H^1(\Omega_0)}
+\|\chi_{\eps}\|_{H^1(\Omega_\eps\setminus\Omega_0)}\eto 0$.

\par\bigskip
\begin{definition}\label{definition-convergence}
We will say that the family of domains $\Omega_\eps$
converges spectrally to $\Omega_0$ as $\eps \to 0$
if the eigenvalues and eigenprojectors
of the Neumann Laplacian behave continuously at $\eps=0$.
That is, for any fixed $n\in \mathbb{N}$
we have that $\la_n^\eps\to \la_n ^0$ as $\eps\to 0$,
and  for each $n\in \mathbb{N}$ such that $\la_n^0<\la_{n+1}^0$
the spectral projections
$P_n^\eps:L^2(\R^N)\to H^1(\Omega_\eps)$,
$P_n^\eps(\psi)=\sum_{i=1}^n(\phi_i^\eps,\psi)_{L^2(\Omega_\eps)}\phi_i^\eps$,
satisfy
$$
  \sup\{\|P_n^\eps(\psi)-P_n^0(\psi)\|_{H^1_\eps},\,\psi\in
  L^2(\R^N),\|\psi\|_{L^2(\R^N)}=1\}\eto 0
  \,.
$$
\end{definition}
\par\bigskip

The convergence of the spectral projections is equivalent to
the following: for each sequence $\eps_k\to 0$ there exists a
subsequence, that we denote again by $\eps_k$ and a complete
system of orthonormal eigenfunctions of the limiting problem
$\{\phi_n^0\}_{n=1}^\infty$ such that
$\|\phi_n^{\eps_k}-\phi_n^0\|_{H^1_{\eps_k}}\to 0$ as $k\to \infty$.

In order to write down the characterization,
we need to consider the following quantity
%
%
\begin{equation}\label{Edefinitiontau}
\tau_\eps=\min_{\stackrel{\phi\in H^1(\Omega_\eps)}{\phi=0 \hbox{
in }\Omega_0}}\frac{\ds\int_{\Omega_\eps}|\nabla
\phi|^2}{\ds\int_{\Omega_\eps}|\phi|^2}
\,.
\end{equation}
Observe that $\tau_\eps$ is the first eigenvalue of the following problem
with a combination of Dirichlet and Neumann boundary conditions:
$$
\left\{
\begin{aligned}
  -\Delta u & = \tau u \,, \quad
  &&\Omega_\eps\setminus \bar \Omega_0 \,,
  \\
  u&=0 \,,
  &&\partial \Omega_0 \,,
  \\
  \frac{\partial u}{\partial n}&=0 \,,
  &&\partial\Omega_\eps \setminus \partial\Omega_0 \,.
\end{aligned}
\right.
$$

We can prove the following,

\begin{proposition}\label{SpectralConvergence}
A necessary and sufficient condition for the spectral convergence
of $\Omega_\eps$ to $\Omega_0$ is
\begin{equation}\label{tau-to-infinity}
\tau_\eps\eto +\infty \,.
\end{equation}
\end{proposition}

\par We refer to \cite{ar2} and \cite{AC} for a proof of this result.

\begin{remark}
The fact that $\Omega_0\subset\Omega_\eps$ can be relaxed. It is enough asking that for
each compact set $K\subset \Omega_0$ there exists $\eps(K)$ such that $K\subset \Omega_\eps$
for $0<\eps\leq \eps(K)$, see \cite{AC}.

\end{remark}

\section{Measure convergence of the domains}\label{Measure-convergence}

In this section we provide an answer to the first question.
Observe that in Proposition \ref{SpectralConvergence}
we do not require that $|\Omega_\eps\setminus\Omega_0|\eto 0$.
However, we have the following

\begin{corollary}
In the situation above if $\Omega_\eps$
converges spectrally to $\Omega_0$,
then necessarily $|\Omega_\eps\setminus\Omega_0|\eto 0$.
\end{corollary}

\par\noindent{\bf Proof.}  This result is proved in \cite{AC} but for the sake of completeness and since it is a simple proof,
we include it in here.

If this were not true then we
will have a positive $\eta>0$ and a sequence ${\eps_k}\to 0$ such
that $|\Omega_{\eps_k}\setminus\Omega_0|\geq \eta$. Let
$\rho=\rho(\eta)$ be a small number such that $|\{x\in
\R^N\setminus\Omega_0,\, \dist(x,\Omega_0)\leq \rho\}|\leq
\eta/2$. This implies that
$|\{x\in \Omega_{\eps_k},\, \dist( x,\Omega_0)\geq \rho\}|\geq \eta/2$.
Let us construct a smooth function $\ga$
with $\ga=0$ in $\Omega_0$,
and $\gamma(x)=1$ for $x\in\R^N\setminus\Omega_0$ with
$\dist(x,\Omega_0)\geq  \rho$. Then obviously $\gamma\in
H^1(\Omega_{\eps_k})$ with $\|\nabla
\ga\|_{L^2(\Omega_{\eps_k})}\leq C$ and
$\|\gamma\|_{L^2(\Omega_{\eps_k})}\geq (\eta/2)^{\frac{1}{2}}$. This
implies that $\tau_{\eps_k}$ is bounded.
Hence it is not true that $\tau_\eps\eto +\infty$
and therefore, from Proposition \ref{SpectralConvergence},
we do not obtain the spectral convergence. \cqd

\par\bigskip
In particular, this result implies that the answer to question {\bf (Q1)}
is affirmative. That is, if we have the convergence
of Neumann eigenvalues and eigenfunctions,
necessarily we have that $|\Omega_\eps\setminus\Omega_0|\eto 0$.

\section{Distance convergence of the domains}\label{distance-convergence}

In this section we will provide an answer to question {\bf (Q2)} and,
as a matter of fact, we will see that the answer is No.
We will prove this by constructing an example of a fixed domain $\Omega_0$
and a sequence
of domains $\Omega_\eps$ with $\Omega_0\subset\Omega_\eps$
with the property that dist$(\Omega_\eps,\Omega_0)$ does not
converges to~$0$,
but the eigenvalues and eigenfunctions of the Laplace operator with Neumann boundary
conditions  in  $\Omega_\eps$ converge to the ones in $\Omega_0$,
see Definition~\ref{definition-convergence}.

As a matter of fact in \cite[Section 5.2]{AC} a very particular example of a
dumbbell domain (two disconnected domains joined by a thin channel) is provided so that
the eigenvalues from the dumbbell converge to the eigenvalues of the two disconnected domains and no spectral contribution
from the channel is observed. In this note we will obtain a family of channels for which the same phenomena occurs, see
Corollary \ref{main2},  and will provide a proof, different from the one given in \cite{AC}.

Let us consider a fixed domain $\Omega_0\subset \R^N$
which satisfies that $\Omega_0\subset \{x\in \R^N, \, x_1<0\}$
and such that
\begin{multline*}
  \Omega_0\cap \{x=(x_1,x')\in \R\times\R^{N-1}, -1<x_1<1,\, |x'|\leq \rho\}
  \\
  = \{x=(x_1,x')\in \R\times\R^{N-1}, -1<x_1<0,\, |x'|\leq \rho\}
\end{multline*}
for some fixed $\rho>0$.

\begin{figure}[htbp] 
   \centering
   \epsfig{file=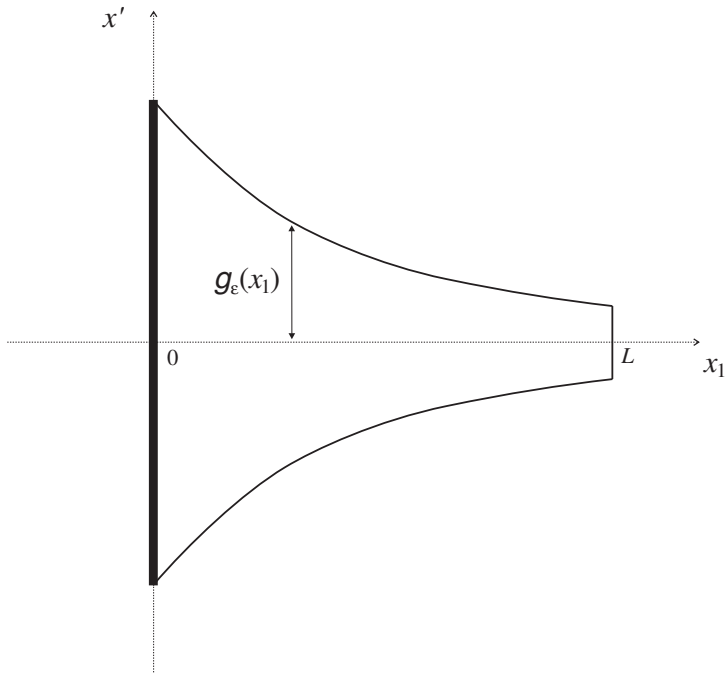,width=0.8\textwidth}
\caption{The exterior perturbation $R_\eps$.
The thick line refers to the supplementary Dirichlet condition
in the problem~\eqref{Equation-for-tau}, while Neumann boundary
conditions are imposed elsewhere.}%
   \label{figure}
\end{figure}

We will construct $\Omega_\eps$ as $\Omega_\eps=$int($\bar\Omega_0\cup\bar R_\eps)$,
where $R_\eps$ is given as follows
\begin{equation}\label{def-R-epsilon}
R_\eps=\{(x_1,x')\in \R\times\R^{N-1}: \, 0<x_1<L, \, |x'|<g_\eps(x_1)\}
\end{equation}
where the function $g_\eps$ will be chosen so that $g_\eps>0$,
$g_\eps\in C^1([0,L])$ and $g_\eps\to 0$
uniformly on $[0,L]$, see Figure~\ref{figure}.
For the sake of notation we denote by
$\Gamma_0^\eps=\partial R_\eps\cap\{ x_1=0\}$
and
$\Gamma_L^\eps=\partial R_\eps\cap\{ x_1=L\}$.

%
We refer to \cite{Rau95} for a general reference
on the behavior of solutions of partial differential equations on thin domains.
See also the recent survey~\cite{Grieser}
for a study on the spectrum of the Laplacian
on thin tubes in various settings,
and for many related references.

Observe that if $L$ is fixed then dist$(\Omega_\eps, \Omega_0)=L$
for each $0<\eps\leq \eps_0$.   Moreover we will show
that for certain choices of $g_\eps$
we obtain the spectral convergence of the Laplace operator.
To prove this results,  we use Proposition \ref{SpectralConvergence}
and show that $\tau_\eps\to +\infty$.
Notice that $\tau_\eps$, defined in (\ref{Edefinitiontau})
is the first eigenvalue of
\begin{equation}\label{Equation-for-tau}
\left\{
\begin{aligned}
  -\Delta u&=\tau u \,,\quad && R_\eps \,,
  \\
  u&=0 \,, \quad && \Gamma_0^\eps \,,
  \\
  \frac{\ds\partial u}{\ds\partial n}&=0 \,,\quad
  &&\partial R_\eps\setminus \Gamma_0^\eps \,.
\end{aligned}
\right.
\end{equation}

Since we have Neumann boundary conditions
on the lateral boundary of~$R_\eps$, there clearly exist
profiles of~$g_\eps$ for which $\tau_\eps$ remains
uniformly bounded as $\eps \to 0$.
In fact, a simple trial-function argument shows that
$\tau_\eps \leq \pi^2/(2L)^2$
whenever $g_\eps(s) \geq g_\eps(0)$ for every $s\in[0,L]$.
The idea to get $\tau_\eps\to +\infty$ consists in choosing
a rapidly decreasing function $s \mapsto g_\eps(s)$,
which enables one to get a large contribution to~$\tau_\eps$
coming from the longitudinal energy
due to the approaching Dirichlet and Neumann boundary conditions
in the limit $\eps \to 0$.
Let us notice that a similar trick to employ the repulsive
contribution of such a combination of the boundary conditions
have been used recently in~\cite{KK} to establish
a Hardy-type inequality in a waveguide;
see also~\cite{K5} for eigenvalue asymptotics
in narrow curved strips with combined Dirichlet and Neumann boundary conditions.
In our case, we are able to show

\begin{proposition}\label{main}
With the notations above, for any function $\gamma\in C^2([0,L])$
satisfying
\begin{equation}\label{gamma-hyp}
0<\alpha_0\leq \gamma\leq \alpha_1<1,
\qquad \dot\gamma(L) \leq 0,
\quad \hbox{ and } \quad
\ddot\gamma \geq \alpha_2>0
\end{equation}
for some positive numbers $\alpha_0$,  $\alpha_1$ and $\alpha_2$,
if we define $g_\eps = \gamma^{1/\eps}$
we have that $\tau_\eps\eto 0$.

In particular, applying Proposition~\ref{SpectralConvergence}
we obtain the convergence of the eigenvalues and eigenfunctions
of the Neumann Laplacian in $\Omega_\eps$ to the ones in $\Omega_0$.
\end{proposition}

\begin{remark}
Observe that a function $\gamma$ satisfying~\eqref{gamma-hyp}
necessarily satisfies that $\dot\gamma(s)<0$ for $0\leq s<L$.
Hence, the function $\gamma$ is decreasing.
\end{remark}

\bproof Since $\tau_\eps$ is given by minimization of the Rayleigh quotient,
$$
\tau_\eps=\inf_{\stackrel{\phi\in H^1(R_\eps)}{\phi=0 \hbox{
in  }\Gamma_0^\eps}}\frac{\ds\int_{R_\eps}|\nabla
\phi|^2}{\ds\int_{R_\eps}|\phi|^2}
$$
we analyze the integral $\int_{R_\eps}|\nabla \phi|^2$
for a smooth real-valued function $\phi$ with $\phi=0$
in a neighborhood of $\Gamma_0^\eps$.
We have
$$
  \int_{R_\eps}|\nabla \phi|^2
  =\int_0^L\int_{|x'|<g_\eps(x_1)}(|\phi_{x_1}|^2+|\nabla_{x'}\phi|^2) \, dx' dx_1
$$

Considering the change of variables  $x_1=y_1$, $x'=g_\eps(y_1)y'$
which transforms $(x_1,x')\in R_\eps$ into $(y_1,y')\in Q$ where
$Q$ is the cylinder $Q=\{(y_1,y'): 0<y_1<L,\, |y'|<1\}$
and performing this change of variables in the integral above,
elementary calculations show that
$$
  \int_{R_\eps}|\nabla \phi|^2
  =\int_Q \left[\left(\varphi_{y_1}
  - \frac{\dot{g}_\eps}{g_\eps} \sum_{i=2}^Ny_i\varphi_{y_i}\right)^2
  +\frac{1}{g_\eps^2}\sum_{i=2}^N|\varphi_{y_i}|^2\right]\, g_\eps^{N-1}dy
$$
where
$\varphi(y)=\phi\big(y_1,g_\epsilon(y_1)y'\big)$.

Writing the above expression in terms of the new function
$\psi(y)=g_\eps(y_1)^\frac{N-1}{2}\varphi(y)$ so that
\begin{align*}
  g_\eps^{(N-1)/2}\varphi_{y_i}
  &=\psi_{y_i}
  \,, \qquad
  i=2,\ldots,N
  \,,
  \\
  g_\eps^{(N-1)/2}\varphi_{y_1}
  &=-\frac{N-1}{2}\frac{\dot{g}_\eps}{g_\eps}\,\psi+\psi_{y_1}
  \,,
\end{align*}
we get,
\begin{eqnarray*}
\lefteqn{
  \int_{R_\eps}|\nabla \phi|^2
  }
  \\
  &&
  =\int_Q \left[\left(-\frac{N-1}{2}\frac{\dot{g}_\eps}{g_\eps}\psi+\psi_{y_1}
  -\frac{\dot{g}_\eps}{g_\eps}\sum_{i=2}^Ny_i\psi_{y_i}\right)^2
  +\frac{1}{g_\eps^2}\sum_{i=2}^N|\psi_{y_i}|^2\right]dy
  \\
  &&
  =\int_Q \Bigg[
  \left(-\frac{N-1}{2}\frac{\dot{g}_\eps}{g_\eps}\psi\right)^2
  + \left(\psi_{y_1}
  - \frac{\dot{g}_\eps}{g_\eps}
  \sum_{i=2}^Ny_i\psi_{y_i}\right)^2
  -(N-1)\frac{\dot{g}_\eps}{g_\eps}\psi\psi_{y_1}
  \\
  &&
  \qquad\qquad
  + (N-1) \frac{\dot{g}_\eps^2}{g_\eps^2}
  \sum_{i=2}^Ny_i\psi_{y_i}\psi
  +\frac{1}{g_\eps^2}\sum_{i=2}^N|\psi_{y_i}|^2\Bigg]dy
  \\
  &&
  \geq \int_Q\left[ \left(\frac{N-1}{2}\right)^2\frac{\dot{g}_\eps^2}{g_\eps^2}\psi^2
  -(N-1)\frac{\dot{g}_\eps}{g_\eps}\psi\psi_{y_1}\right.
  \\
  &&
  \qquad\qquad
  + \left. (N-1)\frac{\dot{g}^2}{g_\eps^2}\sum_{i=2}^Ny_i\psi_{y_i}\psi
  +\frac{1}{g_\eps^2}\sum_{i=2}^N \psi_{y_i}^2
  \right]dy
\end{eqnarray*}
where we have used that
$(\psi_{y_1}-\sum_{i=2}^Ny_i\psi_{y_i}\frac{\dot{g}_\eps}{g_\eps})^2\geq 0$.
Via integration by parts in the second and third term above, we get,
\begin{eqnarray*}
\lefteqn{
  \int_{Q} -(N-1)\frac{\dot{g}_\eps}{g_\eps}\psi\psi_{y_1}dy
  =\int_{|y'|<1}\int_0^L -(N-1)\frac{\dot{g}_\eps}{2g_\eps}(\psi^2)_{y_1}dy_1dy'
}
  \\
  &&
  =\int_{|y'|<1}\left( -\left[(N-1)
  \frac{\dot{g}_\eps}{2g_\eps}\psi^2\right]_{y_1=0}^{y_1=L}
  +\int_0^L(N-1)\left(\frac{\dot{g}_\eps}{2g_\eps}\right)^{\!\!\mbox{\large.}}
  \psi^2dy_1\right)dy'
  \\
  &&
  =-\int_{|y'|<1}(N-1) \frac{\dot{g}_\eps(L)}{2g_\eps(L)}\psi^2(L,y')dy'
  +\int_Q\frac{N-1}{2}\left(\frac{\ddot{g}_\eps}{g_\eps}
  -\frac{\dot{g}_\eps^2}{g_\eps^2}\right)\psi^2dy
\end{eqnarray*}
and
\begin{multline*}
  \int_{Q} (N-1)\frac{\dot{g}_\eps^2}{g_\eps^2}\sum_{i=2}^Ny_i\psi_{y_i}\psi dy
  =\int_0^L(N-1)\frac{\dot{g}_\eps^2}{g_\eps^2}
  \sum_{i=2}^N \int_{|y'|<1}y_i\frac{1}{2}(\psi^2)_{y_i} dy'dy_1
  \\
  =
  \int_0^L\frac{N-1}{2}\frac{\dot{g}_\eps^2}{g_\eps^2}\left( \int_{|y'|=1}\psi^2
  - (N-1)\int_{|y'|<1}\psi^2dy'\right)dy_1
  \,.
\end{multline*}

Hence if we require that $\dot{g}_\eps(L)\leq 0$, we  have,
%
\begin{equation}\label{equation1}
\begin{array}{l}
  \ds \int_{R_\eps}|\nabla \phi|^2
  \geq \int_Q\left[\frac{N-1}{2}\frac{\ddot{g}_\eps}{g_\eps}
  -\left(\left(\frac{N-1}{2}\right)^2+\frac{N-1}{2}\right)
  \frac{\dot{g}_\eps^2}{g_\eps^2}\right] \psi^2dy
  \\
\ds\quad\qquad  +\int_0^L\frac{N-1}{2}\frac{\dot{g}_\eps^2}{g_\eps^2}
  \left( \int_{|y'|=1}\psi^2 dy'\right)dy_1
  +\int_Q \frac{1}{g_\eps^2}\sum_{i=2}^N \psi_{y_i}^2dy
  \,.
\end{array}
\end{equation}
%
The last two terms in this expression can be written as
$$
  \int_0^L\frac{1}{g_\eps^2(y_1)}\left(\int_{|y'|\leq 1}
  |\nabla_{y'}\psi|^2+\frac{N-1}{2} \, \dot{g}_\eps^2(y_1)\int_{|y'|=1}\psi^2\right)dy_1
$$
and we have that
$$
  \int_{|y'|\leq 1} |\nabla_{y'}\psi|^2
  +\frac{N-1}{2}\,\dot{g}_\eps^2\int_{|y'|=1}\psi^2
  \geq \rho \int_{|y'|\leq 1}\psi^2
$$
with $\rho=\rho(y_1)$ being the first eigenvalue of the problem
$$
\left\{
\begin{aligned}
  -\Delta_{y'}\psi&=\rho \psi \,,
  \quad && |y'|<1 \,,
  \\
  \frac{\partial \psi}{\partial n}+\frac{N-1}{2}\,\dot{g}_\eps^2(y_1)\psi
  &=0 \,,
  \quad && |y'|=1 \,,
\end{aligned}
\right.
$$
where~$n$ denotes the outward unit normal vector field
to the $(N-2)$ dimensional unit sphere
$
  S_1=\{y'\in\R^{N-1}:\,|y'|=1\}
$.

We claim that  if we denote by $\lambda(\eta)$ the first eigenvalue of
$$
  \left\{
  \begin{aligned}
  -\Delta_{y'}\psi&=\lambda \psi \,,\quad
  && |y'|<1 \,,
  \\
  \frac{\partial \psi}{\partial n}+\eta \psi&=0 \,,\quad
  && |y'|=1 \,,
  \end{aligned}
  \right.
$$
we have that $\frac{\lambda(\eta)}{\eta}\to \frac{|S_1|}{|B_1|}$ as $\eta\to 0$,
where  $B_1$ is the $(N-1)$ dimensional unit ball and $S_1$ its surface,
which satisfy $|S_1|=(N-1)|B_1|$.
As a matter of fact by standard continuity result
we know that $\lambda(\eta)\to 0$
and its eigenfunction $\psi_\eta$, which is radially symmetric, converges
to the constant function $1/\sqrt{|B_1|}$,
which is the first eigenfunction of the Neumann eigenvalue problem.
But
$$
  \lambda(\eta)
  =\int_{B_1}|\nabla_{y'} \psi_\eta|^2+\eta\int_{S_1}|\psi_\eta|^2
  \geq \eta \int_{S_1}|\psi_\eta|^2
$$
 which implies that
 $$
   \frac{\lambda(\eta)}{\eta}\geq  \int_{S_1}|\psi_\eta|^2\to \frac{|S_1|}{|B_1|}
   \,.
 $$
 Moreover, using $\psi=1/\sqrt{|B_1|}$ as a test function
 in the Rayleigh quotient for $\lambda(\eta)$, we immediately obtain
 $\lambda(\eta)\leq \eta \frac{|S_1|}{|B_1|}$.
 This proves our claim.
 In particular, given $\delta>0$ small, we can choose $\eta_0=\eta_0(\delta)$
 such that $\lambda(\eta)>(N-1-\delta)\eta$ for $0<\eta\leq \eta_0$.

Therefore, if we choose the function $g_\eps$ such that $\dot{g}_\eps(y_1)\to 0$
uniformly in $y_1\in [0,L]$, we have
that $\rho(y_1)\geq \frac{(N-1)(N-1-\delta)}{2}\,\dot{g_\eps}^2(y_1)$
for $\eps$ small enough.

Hence,
\begin{align*}
  \int_{R_\eps}|\nabla \phi|^2
  \ \geq\ &  \int_Q
  \bigg\{\frac{N-1}{2}\frac{\ddot{g}_\eps}{g_\eps}
  -\bigg[
  \left(\frac{N-1}{2}\right)^2
  \\
  & \qquad -\frac{(N-1)(N-1-\delta)}{2}
  +\frac{N-1}{2}
  \bigg] \frac{\dot{g_\eps}^2}{g_\eps^2}
  \bigg\} \, \psi^2dy
  \\
  \ =\ & \frac{N-1}{2}\int_Q
  \left\{\frac{\ddot{g_\eps}}{g_\eps}-\left[\frac{N-1}{2}
  -(N-1-\delta)+1\right]\frac{\dot{g}_\eps^2}{g_\eps^2}
  \right\}\psi^2 dy
\end{align*}
and observe that the number $\kappa=\frac{N-1}{2}-(N-1-\delta)+1$
is strictly less than one for all values of $N\geq 2$
choosing a fixed and small $\delta>0$.
If we denote by
$$
  m_\eps=\inf_{0\leq y_1\leq L}\left(
  \frac{\ddot{g}_\eps}{g_\eps}-\kappa\,\frac{\dot{g}_\eps^2}{g_\eps^2}
  \right)
$$
then
$$
  \int_{R_\eps}|\nabla \phi|^2
  \geq \frac{N-1}{2}\,m_\eps \int_Q\psi^2
  = \frac{N-1}{2}\,m_\eps \int_{R_\eps}\phi^2
  \,.
$$
Consequently, $\tau_\eps \geq \frac{N-1}{2}\,m_\eps$.

Let us see that we can make a choice of the family of functions $g_\eps$,
satisfying the two previous conditions we have imposed,
that is $\dot{g}_\eps(L) \leq 0$ and $\dot{g}_\eps(y_1)\to 0$
uniformly in $0\leq y_1\leq L$
such that $m_\eps\to +\infty$ as $\eps\to 0$.

Let us choose a function $\gamma \in C^2([0,L])$
satisfying (\ref{gamma-hyp})
and let $g_\eps=\gamma^{1/\eps}$.
Then, we have
$$
  \dot{g}_\eps
  = \frac{1}{\eps}\gamma^{\frac{1}{\eps}-1}\dot\gamma
  \,, \qquad
  \ddot{g}_\eps
  = \frac{1}{\eps}(\frac{1}{\eps}-1)\gamma^{\frac{1}{\eps}-2}
  \dot\gamma^2+\frac{1}{\eps}\gamma^{\frac{1}{\eps}-1}\ddot\gamma
  \,,
$$
and simple calculations show that
$$
  \frac{\ddot{g}_\eps}{g_\eps}-\kappa\,\frac{\dot{g}_\eps^2}{g_\eps^2}
  =\left[\frac{1}{\eps}(\frac{1}{\eps}-1)
  -\kappa \left(\frac{1}{\eps}\right)^2\right]
  \left(\frac{\dot\gamma}{\gamma}\right)^2
  +\frac{\ddot\gamma}{\eps\gamma}
  \geq \frac{\alpha_2}{\alpha_0}\frac{1}{\eps}
$$
for $\eps>0$ small enough so that $\frac{1}{\eps}(\frac{1}{\eps}-1)-\kappa \left(\frac{1}{\eps}\right)^2\geq 0$ . This shows that $m_\eps\to +\infty$ and it proves
the proposition. \cqd

\par\bigskip
\begin{remark}
Now that we have been able to construct a thin domain $R_\eps$
as in \eqref{def-R-epsilon} such that $\tau_\eps\eto +\infty$,
we can construct
another thin domain $\tilde R_\eps$ such that its ``length'' goes to infinity,
its width goes to zero and still $\tilde \tau_\eps\eto +\infty$, where
$\tilde\tau_\eps$ is the first eigenvalue of \eqref{Equation-for-tau}
in $\tilde R_\eps$ instead of $R_\eps$.

For this, let $R_\eps$ be a thin domain constructed as in Proposition~\ref{main}
and let $\rho_\eps$ be a sequence with $\rho_\eps\to +\infty$
such that $\frac{\tau_\eps}{\rho_\eps^2}\to +\infty$
and $\alpha_1^{1/\eps} \rho_\eps \to 0$.
Define $\tilde R_\eps=\rho_\eps R_\eps$, that is
$$
  \tilde R_\eps=\{(x_1,x'): \, 0<x_1<\rho_\eps L, \, |x'|<\rho_\eps g_\eps(x_1)\}
  \,,
$$
then $0<\rho_\eps g_\eps(x_1)\leq \alpha_1^{1/\eps} \rho_\eps \eto 0$
and $\tilde\tau_\eps=\frac{\tau_\eps}{\rho_\eps^2}\eto +\infty$.
\end{remark}

Observe that if we require also a Dirichlet
boundary condition in $\Gamma_L^\eps$, we can
relax the conditions on $\gamma$ in Proposition~\ref{main}
and in particular the condition $\dot\gamma(L)\leq 0$ can be dropped.
Hence, we can show,

\begin{corollary}\label{main2}
With the notations above, for any function $\gamma\in C^2([0,L])$
satisfying
\begin{equation*}
0<\alpha_0\leq \gamma\leq \alpha_1<1,
\quad \hbox{ and } \quad
\ddot\gamma \geq \alpha_2>0
\end{equation*}
for some positive numbers $\alpha_0$,  $\alpha_1$ and $\alpha_2$,
if we define $g_\eps = \gamma^{1/\eps}$
we have that $\tilde \tau_\eps\eto 0$, where $\tilde \tau\eps$ is the
first eigenvalue of
\begin{equation*}
\left\{
\begin{aligned}
  -\Delta u&=\tau u \,,\quad && R_\eps \,,
  \\
  u&=0 \,, \quad && \Gamma_0^\eps \cup \Gamma_L^\eps\,,
  \\
  \frac{\ds\partial u}{\ds\partial n}&=0 \,,\quad
  &&\partial R_\eps\setminus (\Gamma_0^\eps \cup \Gamma_L^\eps) \,.
\end{aligned}
\right.
\end{equation*}
\end{corollary}

\bproof
This follows easily by a Neumann bracketing argument.
More precisely, from the hypotheses,
$\dot\gamma$ is a strictly increasing function.
Hence, either $\gamma$ is strictly monotone in $(0,L)$,
or there exists a unique $L^*\in (0,L)$ such that $\dot\gamma(L^*)=0$.

In the first case, if~$\gamma$ is decreasing (respectively increasing)
we substitute the Dirichlet boundary condition at $\Gamma_L^\eps$
(respectively at $\Gamma_0^\eps$) by a Neumann one.
Then the new eigenvalue problem gives rise to~$\tau_\eps$
defined exactly in the same way as~(\ref{Equation-for-tau})
(modulo possibly a mirroring of~$R_\eps$)
and we have $\tilde\tau_\epsilon\geq \tau_\eps\to +\infty$ as $\eps\to 0$.

In the second case,
we cut the domain $R_\eps$ in two domains $R_\eps^0=R_\eps\cap
 \{0<x_1<L^*\}$, $R_\eps^1=R_\eps\cap \{ L^*<x_1<L\}$.  We know that
 $\tilde\tau_\eps\geq \inf\{\tau_\eps^0,\tau_\eps^1\}$,
 where $\tau_\eps^0$ and $\tau_\eps^1$ are the corresponding
 eigenvalues in $R_\eps^0$ and $R_\eps^1$ with a Neumann boundary
 condition imposed at the newly created boundary
 $R_\eps\cap \{x_1=L^*\}$ on both domains.
 In both domains
 we can apply Proposition~\ref{main} as in the first case
 so that $\tau_\eps^0, \tau_\eps^1\eto +\infty$,
 which implies $\tilde \tau_\eps\to 0$. \cqd
%

\begin{remark}
This corollary recovers and generalizes the results from
Section~5.2 in  \cite{AC}.
\end{remark}

\end{document}